\documentclass[12pt,reqno]{amsart}
\usepackage{amssymb}

\setbox0=\hbox{$+$}
\newdimen\plusheight
\plusheight=\ht0
\def\+{\;\lower\plusheight\hbox{$+$}\;}

\setbox0=\hbox{$-$}
\newdimen\minusheight
\minusheight=\ht0
\def\-{\;\lower\minusheight\hbox{$-$}\;}

\setbox0=\hbox{$\cdots$}
\newdimen\cdotsheight
\cdotsheight=\plusheight
\def\cds{\lower\cdotsheight\hbox{$\cdots$}}

\renewcommand{\(}{\left\(}
\renewcommand{\)}{\right\)}
\renewcommand{\[}{\left\[}
\renewcommand{\]}{\right\]}

\renewcommand{\pmod}[1]{\,(\textup{mod}\,#1)}
\numberwithin{equation}{section}
  \theoremstyle{plain}
\newtheorem{theorem}{Theorem}[section]
\newtheorem{lemma}[theorem]{Lemma}

\newtheorem{proposition}[theorem]{Proposition}

\newtheorem{question}[theorem]{Question}

\makeatletter
\def\proof{\@ifnextchar[{\@oproof}{\@nproof}}
\def\@oproof[#1][#2]{\trivlist\item[\hskip\labelsep\textit{#2 Proof of\ #1.}~]\ignorespaces}
\def\@nproof{\trivlist\item[\hskip\labelsep\textit{Proof.}~]\ignorespaces}

\makeatother

\makeatletter
\def\leqalignno#1{\displ@y \tabskip\z@ plus\@ne fil
  \halign to\displaywidth{\hfil$\@lign\displaystyle{##}$\tabskip\z@skip
    &$\@lign\displaystyle{{}##}$\hfil\tabskip\z@ plus\@ne fil
    &\kern-\displaywidth\rlap{$\@lign\hbox{\rm##}$}\tabskip\displaywidth\crcr
    #1\crcr}}
\makeatother

\begin{document}
\title{Ramanujan congruences for a class of eta quotients}
\author{Jonah Sinick}\thanks{The author was partially supported by the NSF}
\maketitle
\begin{abstract}
We consider a class of generating functions analogous to the generating function of the partition function and establish a bound on the primes $\ell$ for which their coefficients $c(n)$ obey congruences of the form $c(\ell n + a) \equiv 0 \pmod \ell$. We apply this result to obtain a complete characterization of the congruences of the same form that the sequences $c_N(n)$ satisfy, where $c_N(n)$ is defined by $ \sum_{n=0}^{\infty} c_N(n)q^n = \prod_{n=1}^{\infty} \frac{1}{(1-q^n)(1-q^{Nn})}$ for $N > 1$. This last result answers a question of H.-C. Chan.
\end{abstract}
\section{Introduction}

Some of Ramanujan's most influential results are his congruences for the partition function $p(n)$ (mod 5), (mod 7) and (mod 11). For $n \geq 1$, the function $p(n)$ is defined to be the number of ways of writing $n$ as a sum of positive integers in non-increasing order. By convention, one sets $p(0) = 1$ and $p(n) = 0$ for $n < 0$. Ramanujan discovered that for any $n \in \mathbb{Z}$, we have
\begin{align}
\begin{cases}
  p(5n+4)& \equiv 0 \pmod{5} \\
  p(7n+5)& \equiv 0 \pmod{7} \\
  p(11n+6)& \equiv 0 \pmod{11}.
\end{cases}
\end{align}
He proved the congruences in (1.1) starting from the fact that $\prod_{n=1}^{\infty} \frac{1}{1-q^n} = \sum_{n=0}^{\infty} p(n)q^n$. The congruences in (1.1) have inspired much research in $q$-series, combinatorics and modular forms. For a short survey of this work, we refer the reader to \cite{ao}.

One noticeable feature of the congruences listed above is that that they all take the form $p(\ell n + a) \equiv 0 \pmod {\ell}$ where $\ell$ is prime. It is natural to ask whether $p(n)$ satisfies any other congruences of the same form. In \cite{ab}, Ahlgren and Boylan showed that Ramanujan's congruences are the only congruences of this form:  if $\ell$ is prime,  $0 \leq a \leq \ell -1$  and  $p(\ell n + a) \equiv 0 \pmod {\ell}$, then $(\ell, a) \in \{(5, 4), (7, 5), (11, 6)\}$.

In \cite{hcc}, H.-C. Chan defined a sequence $r(n)$ by the formula

\begin{align*}
\prod_{n=1}^{\infty} \frac{1}{(1-q^n)(1-q^{2n})} = \sum_{n=0}^{\infty} r(n)q^n
\end{align*}
and proved that $r(3n + 2) \equiv 0 \pmod 3$. The form of this last congruence parallels Ramanujan's three congruences listed above: it is of the form $r(\ell n + a) \equiv 0 \pmod \ell$ for $\ell$ prime. In \cite{hcc2}, Chan asked if there are any other congruences of the same form. In this paper we answer his question in the negative as a consequence of Theorem 1.1 below. Define a \emph{Ramanujan congruence} for a sequence $c(n)$ to be a congruence of the form  $c(\ell n + a) \equiv 0 \pmod \ell$ for all $n \in \mathbb{Z}$ with $\ell$ prime. Without loss of generality we can take $0 \leq a \leq \ell -1$.

\begin{theorem} Let $N  > 1$. Define $c_N(n)$ by
\begin{align*}
\prod_{n=1}^{\infty} \frac{1}{(1-q^n)(1-q^{Nn})} = \sum_{n=0}^{\infty} c_N(n)q^n.
\end{align*}

Let $\ell$ be prime, $0 \leq a \leq \ell - 1$ and suppose that
\begin{align*}
c_N(\ell n + a) \equiv 0 \pmod \ell
\end{align*}
 for all $n$. Then $ 2 < \ell \leq 11$. Moreover,
 \begin{itemize}
\item $\ell = 3$ if and only if $N = 2$ and $a = 2$,
\item $\ell = 5$ if and only if $N \equiv 0 \pmod 5$ and $a = 4$,
\item $\ell = 7$ if and only if $N \equiv 0 \pmod 7$ and $a = 5$,
\item $\ell = 11$ if and only if $N \equiv 0 \pmod {11}$ and $a = 6$.
\end{itemize}
 \end{theorem}

Theorem 1.1 gives a complete characterization of Ramanujan congruences for the family of sequences $c_N(n)$. The reader should note that when $c_N(n)$ satisfies a sufficient condition for the existence of a Ramanujan congruence (mod $\ell$), $\ell = 5, 7$ or $11$, the congruence follows trivially from the known congruences for $p(n)$, so that the effect of Theorem 1.2 is that the sequences $c_N(n)$ obey no Ramanujan congruences other than Chan's and those that come from the Ramanujan congruences for $p(n)$ in a trivial way.

We prove Theorem 1.1  using a more broadly applicable theorem which we now state.

\begin{theorem} Let $S = (a_1, a_2,\ldots,a_{j})$ be a sequence of positive integers with $j$ even and define $c(n)$ by
 \begin{equation}
\prod_{n=1}^{\infty} \prod_{i =1}^{j} \frac{1}{(1-q^{a_{i}n})} = \sum_{n=0}^{\infty} c(n)q^n.
 \end{equation}

Let $N = \operatorname{lcm}(a_1, a_2, \ldots , a_{j})$. Then if $c(n)$ obeys a Ramanujan congruence (mod $\ell$), then $\ell|N$ or $\ell \leq \max(5, j + 4)$.
 \end{theorem}

It follows that if $c(n)$ obeys a Ramanujan congruence (mod $\ell$), then $\ell \leq \max(N, 5, j + 4)$. This finiteness result contrasts with Treneer's result \cite{tre} that there are infinitely many congruences of the form $c(An + B) \equiv 0 \pmod M$ where $A, M \in \mathbb{N}$ are allowed to be arbitrary. Treneer's result is a broad generalization of the celebrated theorem of Ono \cite{onopaper} showing the existence of infinitely many congruences for the partition function $p(n)$ and its extension by Ahlgren \cite{a}. These results are quite a bit sharper than we indicate here; we refer the reader to the original sources for more information.

Upon taking $a_i = 1$ for each $i$, Theorem 1.2 reduces to a result of Kiming and Olsson \cite{ko} that there is an explicit bound on those $\ell$ for which there is a Ramanujan congruence (mod $\ell$) for the coefficients of an even power of the generating function of the partition function. Our method of proof is essentially that of Kiming and Olsson but we do not follow their exposition in detail. Kiming and Olsson used the theory of modular forms (mod $\ell$) for $SL_2(\mathbb{Z}) = \Gamma_1(1)$. To generalize their results we use certain facts about the ring of modular forms (mod $\ell$)  for $\Gamma_1(N)$ which were provided by Gross \cite{gr}.

The upper bound on $\ell$ implied by Theorem 1.2 is very close to being sharp in $j$ and is sharp in $N$: this follows from the unexceptional congruences for even powers of the generating function for $p(n)$ reported on in \cite{ko}, the exceptional Ramanujan congruences (mod $\ell$) for coefficients of odd powers of the generating function of $p(n)$ as reported on in \cite{boy} and a line of elementary algebra to use the latter Ramanujan congruences to produce Ramanujan congruences for $c(n)$ with $j$ even.

The reader may wonder why Theorem 1.2 is stated for even $j$. We suspect that there is an explicit bound on $\ell$ in for odd $j$ as well, however, rigorously establishing an upper bound on $\ell$ for odd $j$ appears to be substantially more difficult than doing so for even $j$. Indeed, even if we take $a_i = 1$ for all $i$, in contrast to the Kiming and Olsson bound on $\ell$ for even $j$, it appears that there is no established bound on $\ell$ for an arbitrary odd $j$ (but see \cite{boy} for substantial partial results on this matter). The results of Sections $3$ and $4$ hold independent of the parity of $j$; these results may be of use in establishing a generalization of Theorem 1.2 that includes the case with $j$ odd.

In Section 2 we state the facts that we need about the ring of modular forms (mod $\ell$) for $\Gamma_1(N)$. In Section 3 we use Lemma 4.1 to determine  determine $a$ if $c(\ell n + a) \equiv 0 \pmod \ell$ and $\ell$ is larger than an explicit bound. In Section 4 we prove Lemma 4.1 which we use in Section 5 to prove Theorem 1.2. In Section 6 we use Theorem 1.2 to prove Theorem 1.1. In Section 7 we conclude with comments and open questions.

\section{Modular Forms (mod $\ell$) for $\Gamma_1(N)$, $N \geq 4$}

Before stating the facts that we need about modular forms (mod $\ell$) for $\Gamma_1(N)$, we define the filtration, the operator $\theta$ and the Eisenstein series for $SL_2(\mathbb{Z})$. A general reference for this material is \cite{ono}.

Given an element $f(z)\in M_k(\Gamma_1(N))\cap\mathbb{Z}[[q]]$ and a prime $\ell \in \mathbb{Z}$, reducing the Fourier expansion of $f(z)$ (mod $\ell$) gives an element $\tilde{f}\in \mathbb{F}_\ell[[q]]$. We call such a series a ``modular form (mod $\ell$) for $\Gamma_{1}(N)$.'' We want a notion of ``weight'' for such a series. At first blush one might attempt to define the weight of such a series as the weight of the preimage under the reduction map, but there are many preimages of any such series and not all have the same weight. This motivates the definition of the filtration of a modular form $f\in M_{k}(\Gamma_1(N))\cap\mathbb{Z}[[q]], f \not \equiv 0 \pmod \ell$
which is defined as follows:
\begin{align*}
w_\ell(f):=\text{min}\{k':\tilde{f}\in\widetilde{M}_{k'}(\Gamma_{1}(N))
\end{align*}
where
\begin{align*}
\widetilde{M}_{k'}(\Gamma_{1}(N))=\{\tilde{f}: f(z)\in M_{k'}(\Gamma_1(N))\cap\mathbb{Z}[[q]]\}.
\end{align*}\
We mildly abuse notation and given $\tilde{f}$ a modular form (mod $\ell$) with $w_\ell(f) = k$, we also call the preimages of $\tilde{f}$ under the reduction map ``modular forms (mod $\ell$) with filtration $k$.''

Given $f(z) = \sum_{n=0}^{\infty} c(n)q^n$ where $q = e^{2\pi i z}$, define
\begin{align*}
\theta{f}:= \frac{1}{2\pi i}\frac{df}{dz} = \sum_{n=0}^{\infty} nc(n)q^n .
\end{align*}
\noindent The Eisenstein series for $SL_2(\mathbb{Z})$ of weight 2k is
\begin{align*}
E_{2k}(z) = 1 - \frac{4k}{B_{2k}}\sum_{n=0}^{\infty} \sigma_{2k-1}(n)q^n
\end{align*}

\noindent For $k > 1$, $E_{2k}(z)$ is a modular form for $SL_2(\mathbb{Z})$ of weight $2k$. For $k = 1$, $E_{2k}(z)$ is not a modular form for $SL_2(\mathbb{Z})$ but rather a quasi-modular form. Given a complex analytic function $f(z)$ defined on the upper half plane and an integer $k > 0$ and
$ \bigl( \begin{smallmatrix}
a&b\\ c&d
\end{smallmatrix} \bigr)
 \in SL_2(\mathbb{Z})$, as usual define the slash operator of weight $k$ by
\begin{align}
f(z)|_{k}M = (cz + d)^{-k}f\left(\frac{az + b}{cz + d}\right).
\end{align}
\noindent Though the slash operator depends on $k$ we often omit the subscript $k$ to avoid cumbersome notation. Returning to our comment about $E_{2}(z)$, as mentioned on pg. 18 of \cite{ds}, if $M$ is as above we have
\begin{equation}
E_2(z)|_{2}M = E_2(z) - \frac{6ic}{\pi(cz + d)}.
\end{equation}
If $f$ is a modular form of weight $k$ for $\Gamma_1(N)$ then $12\theta f - k E_2 f$ is a modular form of weight $k + 2$ for $\Gamma_1(N)$. This is Lemma 3 of \cite{swd} for $N = 1$ and is proved for arbitrary $N$ in exactly the same way as for $N = 1$: by unpackaging the definitions and using (2.2). Theorem 2a) from \cite{swd} is that $E_{\ell - 1} \equiv 1 \pmod \ell$ and $E_{\ell + 1} \equiv E_2 \pmod \ell$. Putting these results together we obtain Lemma 2.1.

\begin{lemma}
If $f\in M_{k}(\Gamma_1(N))\cap\mathbb{Z}[[q]]$, then defining $R$ to be
\begin{equation}
R = \left(\theta f - \frac{k}{12}E_{2}f\right)E_{\ell -1} + \frac{k}{12}E_{\ell + 1}f,
\end{equation}
$R$ is a modular form of weight $k + \ell + 1$ such that $R \equiv \theta f \pmod \ell$. In particular, $\theta f$ is a modular form (mod $\ell$) for $\Gamma_1(N)$. It follows that if $\tilde{f} \not\equiv 0 \pmod \ell$, then $w_\ell(\theta f) \leq w_\ell(f) + \ell + 1$.
\end{lemma}
 With Lemma 2.1 and the preceding setup in mind we cite the remaining facts that we need about modular forms (mod $\ell$) for $\Gamma_1(N)$.
\begin{lemma}
Let $N \geq 4$, let $f, g \in M(\Gamma_1(N))\cap\mathbb{Z}[[q]]$, and let $\ell \geq 5$ be prime. Then

\noindent (i) We have  $w_\ell(\theta f) = w_\ell(f) + \ell + 1$ if and only if $w_\ell(f) \not\equiv 0 \pmod \ell$.\\
(ii) If $f$ and $g$ have weights $k_1$ and $k_2$ respectively and $\tilde{f} \equiv \tilde{g}\not\equiv 0 \pmod \ell$, then
\noindent $k_1 \equiv k_2 \pmod {\ell - 1}$.\\
(iii) If $\ell \nmid N$ then for $i \geq 0$, $ w_\ell(f^i) = i\cdot w_\ell(f).$
\end{lemma}
For a proof of Lemma 2.2, see Section 4 of \cite{gr}.\\

\noindent If $f(z) = \sum_{n=0}^{\infty} c(n)q^n$ where $q = e^{2\pi i z}$, define $f|U_\ell$ by
\begin{align*}
f|U_\ell:= \sum_{n=0}^{\infty} c(\ell n)q^n.
\end{align*}
A crucial elementary fact is that if $f\in M_{k}(\Gamma_1(N))\cap\mathbb{Z}[[q]]$, then there is a relationship between $\theta f$ and $f|U_\ell$:
\begin{align*}
(f|U_\ell)^\ell \equiv f - \theta^{\ell-1}f \pmod \ell.
\end{align*}
\noindent It follows that
\begin{equation}
f|U_\ell \equiv 0 \pmod \ell \iff  \theta^{\ell-1}f \equiv f \pmod \ell.
\end{equation}

\section{Determination of $a$ if $c(\ell n + a) \equiv 0 \pmod \ell$}

Let $\delta_\ell = \frac{\ell^2 - 1}{24}$ and as usual let $\Delta(z) = q\prod_{n=1}^{\infty} (1-q^n)^{24}$.  In this section we prove the following.

\begin{lemma}
Let $c(n)$, $j$, $a_1, ..., a_j$ and $N$ be as in Theorem 1.2, and let $\ell > \max(5, j + 3)$ be a prime such that $\ell\nmid N$.

\noindent (i) Then $c(\ell n + a) \equiv 0 \pmod \ell$ if and only if $d(\ell n + b) \equiv 0 \pmod \ell$, where $d(n)$ is defined by
\begin{align}\sum_{n=0}^{\infty} d(n)q^n = \{\prod_{i =1}^{j}  \Delta(a_{i}z) \}^{\delta_\ell},\end{align}
 and $b$ is defined by $24a \equiv 24b + (\sum_{i=1}^{j} a_i) \pmod \ell$.

\noindent (ii) In part (i) we have $b \equiv 0 \pmod \ell$ so that $24a \equiv (\sum_{i=1 }^{j} a_i) \pmod \ell$.
\end{lemma}

The specific tool that we use is a modified form of Proposition 3 from \cite{ko}. One modification is the addition of an additional hypothesis which is implicitly assumed in the proof of Proposition 3 and not explicitly stated. The other modification is that we replace the space $M_k(\Gamma_1(1))$ in Proposition 3 with $M_k(\Gamma_1(N))$ for $N \geq 4$. This yields a true statement because the proof of Proposition 3 given in \cite{ko} is the same word for word for any $N$ for which Lemma 2.1 and Lemma 2.2 of Section 2 are true.

\begin{proposition}
(After Proposition 3 in \cite{ko}) Let $\ell \geq 5$ be prime and $N \geq 4$, $\ell\nmid N$. Suppose that $f(z) \in M_k(\Gamma_1(N))$ has $\ell$-integral Fourier coefficients, $w_\ell(f(z)) \not\equiv 0 \pmod \ell$, and $\theta(f(z))\not\equiv 0 \pmod \ell$. Suppose further that $w_\ell(\theta^{m} f(z)) \geq w_\ell(f(z))$. Then if the Fourier coefficients $d(n)$ of $f(z)$ satisfy $d(\ell n + b) \equiv 0 \pmod \ell$, one of the following is true: $b = 0$, $w_\ell(f(z))\equiv (\ell + 1)/2 \pmod \ell$ or $w(f(z))\equiv (\ell + 3)/2 \pmod \ell$.
\end{proposition}

The hypothesis that is implicitly assumed in the proof of Proposition 3 of \cite{ko} is that $w_\ell(\theta^{m} f(z)) \geq w_\ell(f(z))$.\\

\noindent\emph{Proof of Lemma 3.1}. Since $\ell > 3$,  $24|(\ell^2 -1)$. Write $-1 = -\ell^2 + (\ell^2 - 1)$. Then we have
\begin{align*}
\sum_{n=0}^{\infty} c(n)q^n &=  \prod_{i =1}^{j} \prod_{n=1}^{\infty} (1-q^{a_{i}n})^{-1} =  \prod_{i =1}^{j} \prod_{n=1}^{\infty} (1-q^{a_{i}n})^{-\ell^2}(1-q^{a_{i}n})^{\ell^2 - 1} \\
&= \prod_{i =1}^{j} \prod_{n=1}^{\infty} \left\{ ( 1-q^{a_{i}n})^{-\ell^2}q^{\frac{-a_{i}(\ell^2-1)}{24}}((q^{a_{i}/24})(1-q^{a_{i}n}))^{\ell^2 - 1} \right\} \\
&= \left\{ \prod_{i =1}^{j}  \prod_{n=1}^{\infty} (1-q^{a_{i}n})\right\}^{-\ell^2} \left\{ \prod_{i =1}^{j} q^{-a_{i} \cdot \delta_\ell} \Delta(a_{i}z)^{\delta_\ell} \right\}. \\ \end{align*}

\noindent It follows that
\begin{equation}
q^{-(\delta_\ell \cdot \sum_{i=1}^{j} a_i)}\left\{\prod_{i =1}^{j}  \Delta(a_{i}z)\right\}^{\delta_\ell} = \left\{\prod_{i =1}^{j} \prod_{n=1}^{\infty} (1-q^{a_{i}n})\right\}^{\ell^2} \left\{ \sum_{n=0}^{\infty} c(n)q^n \right\}.
\end{equation}
\noindent Multiplying (3.1) by $q^{-a}$, applying the operator $U_\ell$ and recalling the definition of $d(n)$ gives
\begin{align*}
\sum_{n=0}^{\infty} d\left(\ell n +  \left(\delta_\ell \cdot \sum_{i=1}^{j} a_i\right)+ a \right)q^n  = \left\{\prod_{i =1}^{j} \prod_{n=1}^{\infty} (1-q^{a_{i}n})\right\}^{\ell} \left\{ \sum_{n=0}^{\infty} c(\ell n + a)q^n \right\}.
\end{align*}
\noindent It follows that
\begin{align*}
 \sum_{n=0}^{\infty} d\left(\ell n +  \left(\delta_\ell \cdot \sum_{i=1}^{j} a_i\right)+ a \right)q^n  \equiv 0 \pmod \ell \iff \sum_{n=0}^{\infty} c(ln + a)q^n \equiv 0 \pmod \ell.
\end{align*}
This completes the proof of Lemma 3.1(i).

In view of Lemma 3.1(i), to prove Lemma 3.1(ii), it suffices to show that if $d(\ell n + b)\equiv 0 \pmod \ell$, then $b \equiv 0 \pmod \ell$. The point is that we can show that $b \equiv 0 \pmod \ell$ using the theory of modular forms (mod $\ell$) for $\Gamma_1(N)$. Indeed, it is a standard fact that $\Delta(a_{i}z)$ is a modular form for $\Gamma_0(a_{i})$ so that $\Delta(a_{i}z)$ is a modular form for $\Gamma_1(a_{i})$ and

\begin{align*}
F_{\ell}(z) = \{\prod_{i =1}^{j}  \Delta(a_{i}z) \}^{\delta_\ell}
\end{align*}

\noindent is a modular form for $\Gamma_1(N)$ where $N =  \operatorname{lcm}(a_1, a_2,\ldots, a_{j})$.

To apply Proposition 3.2, we treat $F_\ell$ as a modular form on $\Gamma_1(N')$ where $N' = N$ if $N > 3$ and $N' = 6$ if $N \leq 3$. We now verify that $F_\ell(z)$ satisfies the hypotheses of Proposition 3.2. Of course $F_\ell(z)$ has $\ell$-integral Fourier coefficients having integer Fourier coefficients. By Lemma 4.1 below, $w_\ell(F_\ell) = j(\ell^2 - 1)/2 \not\equiv 0 \pmod \ell$. Since $F_\ell \not\equiv 0 \pmod \ell$, $\theta{F_\ell} \not \equiv 0 \pmod \ell$. Also by Lemma 4.1, $w_\ell(\theta^{m} F_\ell) \geq w_\ell(F_\ell)$.

Applying Proposition 3.2, we see that if $b \neq 0$, either $w_\ell(F_\ell) = j(\ell^2 - 1)/2 \equiv \frac{\ell + 1}{2} \pmod \ell$ or $w_\ell(F_\ell) = j(\ell^2 - 1)/2 \equiv \frac{\ell + 3}{2} \pmod \ell$, but neither possibility occurs since $\ell > j + 3$ by hypothesis. So $b = 0$ as claimed. \qedsymbol

\section{A lemma about $\Theta^m F_\ell$ (mod $\ell$)}

In this section we prove a lemma which we used in the proof of Lemma 3.1 and which we will use further in the proof of Theorem 1.2.

\begin{lemma}
If $m \geq 1$ and $\ell > 3$ is a prime, then
\begin{align*}
w_\ell(\theta^{m} F_\ell) \geq w_\ell(F_\ell) = \frac{j(\ell^2 - 1)}{2}.
\end{align*}
\end{lemma}

Lemma 4.1 appears in \cite{ko} for $N = 1$. The situation is more subtle for a general $N$ than it is for $N = 1$. While to prove Lemma 4.1 for $N = 1$ suffices to consider the Fourier expansion of $F_\ell$ (mod $\ell$) at $\infty$, for a general $N$, $\Gamma_1(N)$ has multiple cusps and we find it necessary to consider the Fourier expansions of $F_\ell$ at each cusp of $\Gamma_1(N)$.

Enumerate the cosets of $\Gamma_1(N)$ in $SL_2(\mathbb Z)$ with $\{i\}_{i = 1, ... 2d_N}$. Let $M_i$ be a representative of the $i$'th coset. Let $\alpha_i$ be the cusp that $M_i$ sends to $\infty$. Denote the minimal period of $F_\ell|M_i$ by $t_i$. Then $F_\ell|M_i$ has a Fourier expansion in powers of $q_{t_i} = e^{\frac{2\pi i z}{t_i}}$. The order of vanishing of $F_\ell$ at $\alpha_i$ is then defined to be the index of the first nonvanishing Fourier coefficient of $F_\ell$ in powers of $q_{t_i}$ and is denoted $\text{ord}_{\alpha_i}(f(z))$.

The Fourier expansions of $F_\ell$ about cusps other than $\infty$ need not have coefficients in $\mathbb{Z}$, but by the $q$-expansion principle, for $N' > 4$ if the Fourier expansion of a modular form $f$ for $\Gamma_1(N')$ about $\infty$ has integer coefficients, then the Fourier coefficients of $f$ about another cusp must lie in $\mathbb{Q}(\zeta_N)$ where $\zeta_N$ is a primitive $N$'th root of unity and have uniformly bounded denominators (see section 12.3 of \cite{di}). This fact is known as the ``bounded denominator property.'' To apply the bounded denominator property we view $F_\ell$ as a modular form for $\Gamma_1(N')$ where $N' = N$ if $N > 4$ and $N' = 12$ if $N \leq 4$.

Before proceeding, we make a remark about the first few paragraphs of Section 2. Rather than considering an element $f$ of $M_{k}(\Gamma_1(N))\cap\mathbb{Z}[[q]]$ and reducing $f$ (mod $\ell$) for some rational prime $\ell$ we can consider elements $g$  of $M_{k}(\Gamma_1(N))\cap L[[q]]$ where $L$ is an algebraic number field and reduce $g$ (mod $v$) for any prime $v$ such that the $v$-adic valuation  of $g$ is $0$. This defines the notion of a ``modular form (mod $v$)'' and allows us to define the filtration $w_v$ for nonvanishing modular forms (mod $v$) in the obvious way. Define the $v$-adic valuation of a power series with coefficients in $L$ to be the minimum of the $v$-adic valuations of the coefficients of the power series (this minimum exists by the bounded denominator property). If we modify the statement of Lemma 2.1 by replacing $\mathbb{Z}[[q]]$ by $L[[q]]$ and replace the modulus of reduction by $v$ where $v$ is a prime above $\ell$ such that the $v$-adic valuation of $f$ is $0$, then the modified Lemma 2.1 is true. We use these facts with $L = \mathbb{Q}(\zeta_N)$. Define $\widetilde{\text{ord}}_{\alpha_i}(f(z))$ to be the order of vanishing of $f$ (mod $v$) at the cusp $\alpha_i$.

As a preliminary to the proof of Lemma 4.1 we prove the following.

\begin{lemma}
Let $m \geq 1$ be an integer and let $v$ be a prime in $\mathbb{Z}[\zeta_N]$ such that $v \nmid 2, 3, N$. Let $f(z)$ be a modular form for $\Gamma_1(N)$ such that $f(z)|M_i$ has coefficients in $\mathbb{Q}(\zeta_N)$ and $v$-adic valuation 0. Let $\alpha_i$ be a cusp of $\Gamma_1(N)$. Then
\begin{align*}
\widetilde{\text{ord}}_{\alpha_i}(\theta^{m}f) \geq \widetilde{\text{ord}}_{\alpha_i}(f).
\end{align*}
\end{lemma}

\noindent\emph{Proof of Lemma 4.2}. By induction it suffices to prove the claim for $m = 1$. Since $R$ in (2.2) satisfies $R \equiv \theta f \pmod v$ for $v \nmid 2,3$ it suffices to show that $\widetilde{\text{ord}}_{\alpha_i}(R) \geq \widetilde{\text{ord}}_{\alpha_i}(f)$. Take $M_i$ to be as in the discussion preceding the statement of Lemma 4.1. Let $k$ be the weight of $f(z)$. Applying the slash operator $|M_i$ of weight $k + \ell + 1$ to both sides of (2.3), we obtain
\begin{align}
R|M_i & = \left(\left(\theta f - \frac{k}{12}E_{2}f \right)E_{\ell -1}\right)|M_i + \left(\frac{k}{12}E_{\ell + 1}f\right)|M_i \notag\\
& = \left((\theta{f})|M_i - \frac{k}{12}(E_2|M_i)(f|M_i)\right)(E_{\ell -1}|M_i) + \frac{k}{12}\left(E_{\ell + 1}|M_i\right)\left(f|M_i\right).
\end{align}

In the second line of (4.1) and in what follows, the slash operators applied to $\theta f$, $E_2$, $f$, $E_{\ell -1}$ and $E_{\ell + 1}$ are of weights $k + 2$, $2$, $k$, $\ell - 1$, and $\ell + 1$ respectively. Now since $E_{\ell -1}$ and $E_{\ell + 1}$ are modular forms for $SL_2(\mathbb Z)$, equation (4.1) becomes
\begin{align}
R|M_i = \left((\theta{f})|M_i - \frac{k}{12}(E_2|M_i)(f|M_i)\right)E_{\ell -1} + \frac{k}{12}\left(E_{\ell + 1}\right)\left(f|M_i\right).
\end{align}
Next we find an alternate expression for $(\theta{f})|M_i$. Applying $\theta$ to both sides of the equation (2.1) gives
\begin{align}
(\theta{f})|M_i = \theta(f|M_i) + \frac{kc}{2\pi i}(cz + d)^{-1}(f|M_i).
\end{align}
Replacing $(\theta{f})|M_i $ in (4.2) by the righthand side of (4.3) and replacing $E_2|M_i$ in (4.2) by the righthand side of (2.2), after simplification (4.2) becomes
\begin{align*}
R|M_i = \left(\theta({f}|M_i) - \frac{k}{12}E_{2}\cdot(f|M_i)\right)E_{\ell -1}  + \frac{k}{12}(E_{\ell + 1})(f|M_i).
\end{align*}
Since $v \nmid N$ and $f|M_i$ has $v$-adic valuation $0$, the Fourier expansion of $\theta({f}|M_i)$ has $v$-adic valuation $0$. It is clear from the definition of $\theta$ that the index of the Fourier first coefficient of $\theta(f|M_i)$ that is nonvanishing (mod $v$) is no smaller than the first Fourier coefficient of $f|M_i$ that is nonvanishing (mod $v$). But then the index of the first Fourier coefficient of $R|M_i$ that is nonvanishing (mod $v$) is no smaller than the first Fourier coefficient of $f|M_i$ that is nonvanishing (mod $v$). This completes the proof. \qedsymbol

\noindent\emph{Proof of Lemma 4.1}. First we prove that $w(F_\ell) = \frac{j(\ell^2 - 1)}{2}$. Consider the functions  $F_\ell|M_i$ for $i \in \{1, ..., 2d_N\}$. Since the $v$-adic valuation of each $F_\ell|M_i$ is finite, for each $i$ there exists a $\beta_i \in \mathbb{Q}$ such that $(\beta_i F_\ell)|M_i$ has $v$-adic valuation 0. Now consider

\begin{align*}
G(z) := \prod_{i=1}^{2d_N}(\beta_{i}F_\ell)|M_i.
\end{align*}

\noindent Since $F_\ell$ is a modular form of weight $\frac{j(\ell^2 - 1)}{2}$ for $\Gamma_1(N)$ and the $M_i$'s are a complete set of representatives of cosets of $\Gamma_1(N)$ in $SL_2(\mathbb{Z})$, $G(z)$ is a modular form of weight $d_{N}j(\ell^2 - 1)$ for $SL_2(\mathbb{Z})$. Let $v$ be a prime above $\ell$ in $\mathbb{Z}[\zeta_N]$.

Since $F_\ell$ is zero-free on $\mathbb{H}$,  $F_\ell|M_i$ is zero-free on  $\mathbb{H}$, so $G(z)$ is zero-free on  $\mathbb{H}$ and the zeros of $G(z)$ all occur at $\infty$. So $G(z)$ must be a nonzero constant multiple of $\Delta(z)^{e}$ where $e = \frac{d_{N}j(\ell^2 - 1)}{12}$. Moreover, by our choice of $\beta_i$, this constant must be nonvanishing (mod $v$). So $w_v(G(z)) = 12e$. But then $w_v(F_\ell) = \frac{j(\ell^2 - 1)}{2}$ and since $F_\ell \in \mathbb{Z}[[q]]$, $w_\ell(F_\ell) = \frac{j(\ell^2 - 1)}{2}$. So we need only show that $w_\ell(\theta^{m} F_\ell) \geq \frac{j(\ell^2 - 1)}{2}$. Since $G$ is a nontrivial multiple of $\Delta^{e}$, it must be that $\widetilde{\text{ord}}_{\infty}(G) = e$. Since $\widetilde{\text{ord}}_{\alpha_i}$ is defined in terms of powers of $q_{t_i}$ while $\widetilde{\text{ord}}_{\infty}$ is defined in terms of powers of $q_{t_1}$ we have $\widetilde{\text{ord}}_{\infty}(\beta_i F_\ell|M_i) = \frac{\widetilde{\text{ord}}_{\alpha_i}(\beta_i F_\ell)}{t_i}$. From the definition of $G$ we see that $\sum_{i=1}^{2d_N} \frac{\widetilde{\text{ord}}_{\alpha_i}(\beta_iF_\ell)}{t_{i}} = \widetilde{\text{ord}}_{\infty}(G(z))$. So
\begin{align*}
\sum_{i=1}^{2d_N} \frac{\widetilde{\text{ord}}_{\alpha_i}(\beta_iF_\ell)}{t_{i}} = e.
\end{align*}
\noindent Now consider
\begin{align*}
H = \prod_{i=1}^{2d_N} (\beta_i \theta^m F_\ell)|M_i.
\end{align*}
\noindent Then $H$ is a modular form (mod $v$) for $SL_2(\mathbb Z)$. We have
\begin{align*}
\widetilde{\text{ord}}_{\infty}(H) = \sum_{i=1}^{2d_N} \frac{\widetilde{\text{ord}}_{\alpha_i}(\beta_{i}\theta^m F_\ell)}{t_i} \geq \sum_{i=1}^{2d_N} \frac{\widetilde{\text{ord}}_{\alpha_i}(\beta_{i}F_\ell)}{t_i} = e,
\end{align*}
where the inequality is a consequence of Lemma 4.2.

The definition of the $\beta_i$ forces $(\beta_i\theta^{m} F_\ell)|M_i \not \equiv 0 \pmod v$, so $H \not\equiv 0 \pmod v$. By Sturm's theorem \cite{st}, $w_v(H) \geq 12e$. But then we see that $w_v(\theta^{m}F_\ell) \geq \frac{j(\ell^2 - 1)}{2}$. Since $\theta^{m}F_\ell \in \mathbb{Z}[[q]]$, we deduce that $w_\ell(\theta^{m}F_\ell) \geq  \frac{j(\ell^2 - 1)}{2}$, completing the proof. \qedsymbol

\section{Proof of Theorem 1.2}

Let $\ell > \max(5, j + 4)$, $\ell \nmid N$. By Lemma 3.1, to prove Theorem 1.2 it suffices to show that $d(\ell n)\equiv 0 \pmod \ell$ leads to a contradiction. Note that it follows from (2.3) that $d(\ell n)\equiv 0 \pmod \ell$ implies $\theta^{\ell-1}F_\ell \equiv F_\ell \pmod \ell$. We analyze the consequences that this has for the sequence $w_\ell(\theta^i f)$, $i \in \{1, 2, \ldots, \ell-1\}$. We will see that the congruence $\theta^{\ell-1}f \equiv f \pmod \ell$ leads to the existence of an $m$ violating the conclusion of Lemma 4.1

To proceed, we need information about the possible sequences $w_\ell(\theta^i F_\ell)$,\\
$i \in \{1, 2, \dots, \ell -1\}$ in terms of $w_\ell(F_\ell)$. Proposition 2 of \cite{ko} is false as stated in \cite{ko} (see appendix) but is true when an additional hypothesis is added to the statement of Proposition 2 and the proof in \cite{ko} is valid once the hypothesis is added. The proof of this modified Proposition 2 carries through without modification when $SL_2(\mathbb{Z})$ is replaced by $\Gamma_1(N)$ for any $N$ for which Lemma 2.1 and Lemma 2.2 hold. For context recall Lemma 2.1 and Lemma 2.2.

\begin{proposition}
(After Proposition 2 of \cite{ko}) Let $\ell \geq 5$ be prime and $N \geq 4$, $\ell\nmid N$. Suppose that $f(z) \in M_k(\Gamma_1(N))$ has $\ell$-integral Fourier coefficients, $w_\ell(f(z)) \not\equiv 0 \pmod \ell$ and $\theta(f(z))\not\equiv 0 \pmod \ell$. Suppose further that $w_\ell(\theta^{m} f(z)) \geq w_\ell(f(z))$. Let $i_1 < i_2 < \ldots < i_v$ be those $i$ with $0 \leq i \leq \ell -1$ for which $w_\ell(\theta^{i}f) \equiv 0 \pmod l$. Write $w_\ell(\theta^{i_j + 1}f) = w_\ell(\theta^{i_j}f) + (\ell + 1) - s_{j}(\ell - 1)$. Write $k = w_\ell(f)$ and let $k_0 \in \{1, 2, \ldots, \ell - 1\}$ be such that $k \equiv -k_0 \pmod \ell$. Then one of the four cases below holds:
 \begin{itemize}
\item (I) $k \equiv 1 \pmod \ell$, $v = 1$, $i_1 = \ell - 1$, and $s_1 = \ell + 1$
\item (II) $k \equiv 2 \pmod \ell$, $v = 1$, $i_1 = \ell - 2$, and $s_1 = \ell + 1$
\item (III) $k \not \equiv 1 \pmod \ell$, $v = 2$, $(i_1,i_2) = (k_0, \ell - 1)$, and $(s_1,s_2) = (k_0 + 1, \ell - k_0)$
\item (IV) $k \not \equiv 1 \pmod \ell$, $v = 2$, $(i_1,i_2) = (k_0, \ell - 2)$, and $(s_1,s_2) = (k_0 + 2, \ell - k_0 - 1)$
\end{itemize}
We have $w_\ell(f) = w_\ell(\theta^{\ell - 1}f)$ if and only if case (II) or case (IV) holds.
\end{proposition}

The necessary hypothesis that is missing in the statement of Proposition 2  is that  $w_\ell(\theta^{m} f) \geq w_\ell(f)$. For a counterexample to the original statement, let $f(z) = \Delta(z)$ and take $\ell = 5$. By Lemma 2.1, there exists a modular form $g(z) \in M_{18}(\Gamma_1(N))$ such that $g \equiv f \pmod 5$. By the equality in Lemma 4.1, $w_5(f) = 12$. Applying Lemma 2.2(i) then gives $w_5(\theta f) = 18$, $w_5(\theta^2 f) = 24$ and $w_5(\theta^3 f) = 30$. Applying $w_5$ to both sides of $\theta^5 f \equiv \theta f \pmod \ell$ forces $w_5(\theta^4 f) = 12$, so that $v = 1$ for $f(z)$ which implies $v = 1$ for $g$. The function $g$ satisfies the hypotheses of the original Proposition 2, but  $w_5(g) = 18 \equiv 3 \pmod 5$ so that if the conclusion of Proposition 2 is true for $g(z)$ then $v = 2$ for $g(z)$, and as we just saw this is not the case.

Now we prove Theorem 1.2. We verified that $F_\ell$ satisfies the hypotheses of Proposition 5.1 in our proof of Lemma 3.1. Taking $f = F_\ell$ in Proposition 5.1 we see that $\theta^{\ell-1}F_\ell \equiv F_\ell \pmod \ell$ implies that we are in case (II) or case (IV) of Proposition 5.1. Actually, we cannot be in case (II) of Proposition 5.1 since Lemma 4.1 shows that  $w_\ell(F_\ell) = j(\ell^2 - 1)/2$ and $\ell > j + 4$, so we are in case (IV) of Proposition 5.1. This implies that if we take $k_0 \equiv -(j)(\ell^2 - 1)/2 \pmod l$, then $w_\ell(\theta^{k_0 + 1} F_\ell) = w_\ell(F_\ell) + (\ell + 1)(k_0 + 1) - (k_0 + 2)(\ell -1) = w(F_\ell) + 2k_0 + 3 - \ell$. We can determine $k_0$ as follows: we have $2k_0 \equiv j \pmod \ell$ so since $j$ is \textbf{even} and $\ell > j$, it must be that  $k_0 = j/2$. So $w_\ell(\theta^{k_0 + 1} F_\ell) = w(F_\ell) + j + 3 - \ell$ and since $\ell > j + 3$, we have $w(\theta^{k_0 + 1} F_\ell) < w(F_\ell)$, contradicting Lemma 4.1 and proving Theorem 1.2. \qedsymbol

\section{Proof of Theorem 1.1}

Let $c_N(n)$ be as in the statement of Theorem 1.1, and assume that $c_N(\ell n + a) \equiv 0 \pmod \ell$. Then by Theorem 1.2, we may assume that $\ell \leq 5$ or $\ell \mid N$. First suppose that $\ell|N$. Write
\begin{align*}
\sum_{n=0}^{\infty} p(n)q^n = \prod_{n=1}^{\infty} \frac{1}{1-q^n} = \left(\sum_{n=0}^{\infty} c_N(n)q^n\right)\left(\prod_{n=1}^{\infty}(1 - q^{Nn})\right).
\end{align*}
Since $\ell|N$ we can write $\prod_{n=1}^{\infty}(1 - q^{Nn}) = \sum_{n=0}^{\infty}y(n)q^{\ell n}$ so that (5.1) becomes
\begin{equation}
\sum_{n=0}^{\infty} p(n)q^n = \left(\sum_{n=0}^{\infty} c_N(n)q^n\right)\left(\sum_{n=0}^{\infty}y(n)q^{\ell n}\right).
\end{equation}
Multiplying (5.1) by $q^{-a}$ and applying $U_\ell$ to both sides gives
\begin{align*}
\sum_{n=0}^{\infty} p(\ell n + a)q^n = \left(\sum_{n=0}^{\infty} c_N(\ell n + a)q^n\right)\left(\sum_{n=0}^{\infty}y(n)q^{n}\right).
\end{align*}
Since $y(0) = 1$, we have $c_N(\ell n + a) \equiv 0 \pmod \ell$ if and only if $p(\ell n + a) \equiv 0 \pmod \ell$ from which it follows that $(\ell, a) \in \{(5, 4), (7, 5), (11, 6)\}$ by the result from \cite{ab} quoted in the Section 1. This establishes Theorem 1.1 assuming that $\ell|N$. So we need only establish Theorem 1.1 assuming that $\ell \leq 5$.

If $\ell \leq 5$ then since $c_N(n) = p(n)$ for $n \leq N$ a short computation shows unless $N$ is as in the bulleted portion of the conclusion of Theorem 1.1, $N \leq 5$. Another short computation together with Chan's result for $N = 2$ show that Theorem 1.1 holds for $\ell \leq 5$.
\qedsymbol

\section{Conclusion}

\noindent In light of our results it is natural to ask:

\begin{question}
Let $c(n)$ be given by $\prod_{n=1}^{\infty} \prod_{i =1}^{j} \frac{1}{(1-q^{a_{i}n})} = \sum_{n=0}^{\infty} c(n)q^n$ where  $j$ is \textbf{odd}. Are there only finitely many $\ell$ for which there is a Ramanujan congruence (mod $\ell$) for $c(n)$? Can one give an explicit bound on $\ell$ if this is so?
\end{question}

\noindent In \cite{boy}, Boylan treated many cases where $j$ is odd and $a_i = 1$ for all $i$. Boylan also reported on the existence of several infinite families of pairs ($j$, $\ell$) such that the coefficients of $\prod_{n=1}^{\infty} \frac{1}{(1-q^n)^j}$ obey a Ramanujan congruence (mod $\ell$), but remarks that there are some pairs ($j$, $\ell$) that do not fit into these families for which there is nevertheless a Ramanujan congruence. A complete characterization of the pairs ($j$, $\ell$) for which there is a Ramanujan congruence appears to be absent from the literature. So we ask the following:

\begin{question}

 Can one give a complete characterization of all tuples ($\ell; a_1, a_2, \ldots, a_j$) for which $c(n)$ given by $\prod_{n=1}^{\infty} \prod_{i =1}^{j} \frac{1}{(1-q^{a_{i}n})} = \sum_{n=0}^{\infty} c(n)q^n$ obeys a Ramanujan congruence (mod $\ell$)?
\end{question}

\noindent While it seems likely that the answer to both parts of Question 7.1 can be answered in the affirmative, the extent of the phenomenon of there being only finitely Ramanujan congruences for the Fourier coefficients of a modular form is quite unclear, motivating:

\begin{question}
Is there a characterization of those weakly holomorphic modular forms $f(z)$ for congruence subgroups of  $SL_2(\mathbb{Z})$ with integer Fourier coefficients such that the Fourier coefficients of $f(z)$ obey only finitely many Ramanujan congruences?
\end{question}

\section{Acknowledgements}

The author was partially supported by the US NSF grant DMS-0707136. I would like to thank Jeremy Rouse for suggesting the problem that led to this paper and for his guidance during this project. In particular, I thank him for substantial help with Section 4. I would like to thank Atul Dixit for being an inspiring study partner during our reading course on modular forms. I thank Michael Dewar for a very careful reading of this paper and for pointing out that Proposition 2 of \cite{ko} is false as stated but easily modified to be applicable in our context. Finally, I thank the anonymous referee for his or her suggestions.\\

\vspace{10 mm}

Department of Mathematics\

University of Illinois at Urbana-Champaign\

1409 W. Green Street, Urbana, IL 61801\\

\emph{Email Address}: jsinick2@math.uiuc.edu

\end{document}